\documentclass[12pt,a4paper,leqno]{article}

\usepackage{amssymb,exscale}
\usepackage[centertags]{amsmath}
\usepackage{graphicx}
\usepackage{slashbox}

\usepackage[dvips]{hyperref}

\renewcommand{\(}{\bigl(}
\renewcommand{\)}{\bigr)\vphantom{)}}

\newcommand{\eps}{\varepsilon}
\newcommand{\ga}{\gamma}
\newcommand{\Om}{\Omega}
\newcommand{\al}{\alpha}
\newcommand{\F}{\mathcal F}

\newcommand{\R}{\mathbb R}

\newcommand{\sif}{$\sigma$\nobreakdash-field}

\newcommand{\SLE}{$ \rm{SLE}_6 $}

\begin{document}

\title{Percolation, boundary, noise: an experiment} 

\author{Boris Tsirelson}

\date{}
\maketitle

\begin{abstract}
The scaling limit of the critical percolation, is it a black noise? The answer
depends on stability to perturbations concentrated along a line. This text,
containing no proofs, reports experimental results that suggest the
affirmative answer.
\end{abstract}

\section*{Introduction}
By \emph{percolation} I mean the two-dimensional critical site percolation on
the triangular lattice (see for instance \cite{La}). Its (full) scaling limit
exists, is conformally invariant and may be described by a random countable
set of non-crossing \SLE\ curves (see \cite{CM}).

By a \emph{noise} I mean a homogeneous continuous product of probability
spaces, as defined by \cite[Def.~3d1]{Ts04}. The scaling limit of percolation
could lead to a noise (and the noise should be black, as defined by
\cite[Def.~7a1]{Ts04}); this idea was discussed informally \cite[Question
11b1]{Ts04} but encountered the following difficulty.

Let a smooth domain $ D \subset \R^2 $ be split by a smooth curve in two
domains $ D_1, D_2 $. The random countable system of non-crossing \SLE\ curves
in $ D $, --- call it the random \emph{configuration} in $ D $, --- being
restricted to $ D_1 $ and $ D_2 $ gives two independent random configurations
in $ D_1, D_2 $. Sewing them together appears to be a subtle matter. Does it
admit any freedom? In other words, is the configuration in $ D $ uniquely
determined by its restrictions to $ D_1 $ and $ D_2 \, $?

More formally: the limiting model on the whole plane $ \R^2 $ is described by
a probability space $ (\Om,\F,P) $ (of configurations). For every open
interval $ (s,t) \subset \R $ we consider the domain $ (s,t) \times \R \subset
\R^2 $ (strip) and the corresponding `local' sub-\sif\ $ \F_{s,t} \subset \F
$. Clearly, $ \F_{r,s} $ and $ \F_{s,t} $ are independent for $ r < s < t $,
and the \emph{upward continuity} is satisfied:
\[
\F_{s,t} = \bigvee_{\eps>0} \F_{s+\eps,t-\eps}
\]
(that is, $ \F_{s,t} $ is the least sub-\sif\ containing all $
\F_{s+\eps,t-\eps} $). The question is, whether
\begin{equation}\label{eq1}
\F_{r,t} = \F_{r,s} \vee \F_{s,t}
\end{equation}
or not; that is, whether the least sub-\sif\ $ \F_{r,s} \vee \F_{s,t} $
containing $ \F_{r,s} $ and $ \F_{s,t} $ is the whole $ \F_{r,t} $, or its
proper sub-\sif. The latter case would mean that the line $ \{s\} \times \R
\subset \R^2 $ cannot be neglected. Note that triviality of the \sif\ $
\F_{s-,s+} = \cap_{\eps>0} \F_{s-\eps,s+\eps} $ is only necessary, since in
general $ \cap_{\eps>0} \( \F_{r,s} \vee \F_{s-\eps,s+\eps} \vee \F_{s,t} \) $
may exceed $ \F_{r,s} \vee \( \cap_{\eps>0} \F_{s-\eps,s+\eps} \) \vee
\F_{s,t} $ (see \cite{Wei} for detailed analysis). Note also that the
difficulty  cannot be avoided by waiving the upward continuity, since the
latter holds for every noise \cite[Prop.~3d3]{Ts04}. We may reformulate
\eqref{eq1} as
\begin{equation}\label{eq2}
\F_{r,t} = \bigvee_{\eps>0} \( \F_{r,s-\eps} \vee \F_{s+\eps,t} \) \, .
\end{equation}
The question is, whether percolation (or rather, its scaling limit) is stable
under a strong perturbation in an infinitesimal strip, or not. This strong
concentrated perturbation is quite different from the distributed weak
perturbation examined in \cite{BKS}. Sensitivity of percolation established
there is \emph{micro-sensitivity} in terms of \cite[Sect.~5c/5.3]{Ts03}, where
it is opposed to \emph{block sensitivity.} The latter also corresponds to a
distributed weak perturbation, but this time the correlation length of the
perturbation tends to $ 0 $ slower than the pitch of the lattice. Block
sensitivity of percolation is noted in \cite[8a2/Remark 8.2]{Ts03}.

Two possibilities remain open. Maybe, percolation is \emph{strip stable,} that
is, satisfies \eqref{eq2}. Then it leads to a noise, and the noise is black
due to the block sensitivity. Or maybe percolation is strip unstable (that is,
violates \eqref{eq2}). Then it does not lead to any noise. A trivial example
of strip instability (and moreover, strip sensitivity) is given by the worst
Boolean function, --- the product of random signs (over all lattice points
in the given domain). Here, trying scaling limit, we face an obstacle, --- a
continuum of independent random signs.

\section[]{\raggedright The idea of the experiment}
\label{sect1}
\begin{figure}[h]
\setlength{\unitlength}{1.7cm}
\begin{picture}(7,3.5)
\put(0.24,0.53){\includegraphics[scale=0.15]{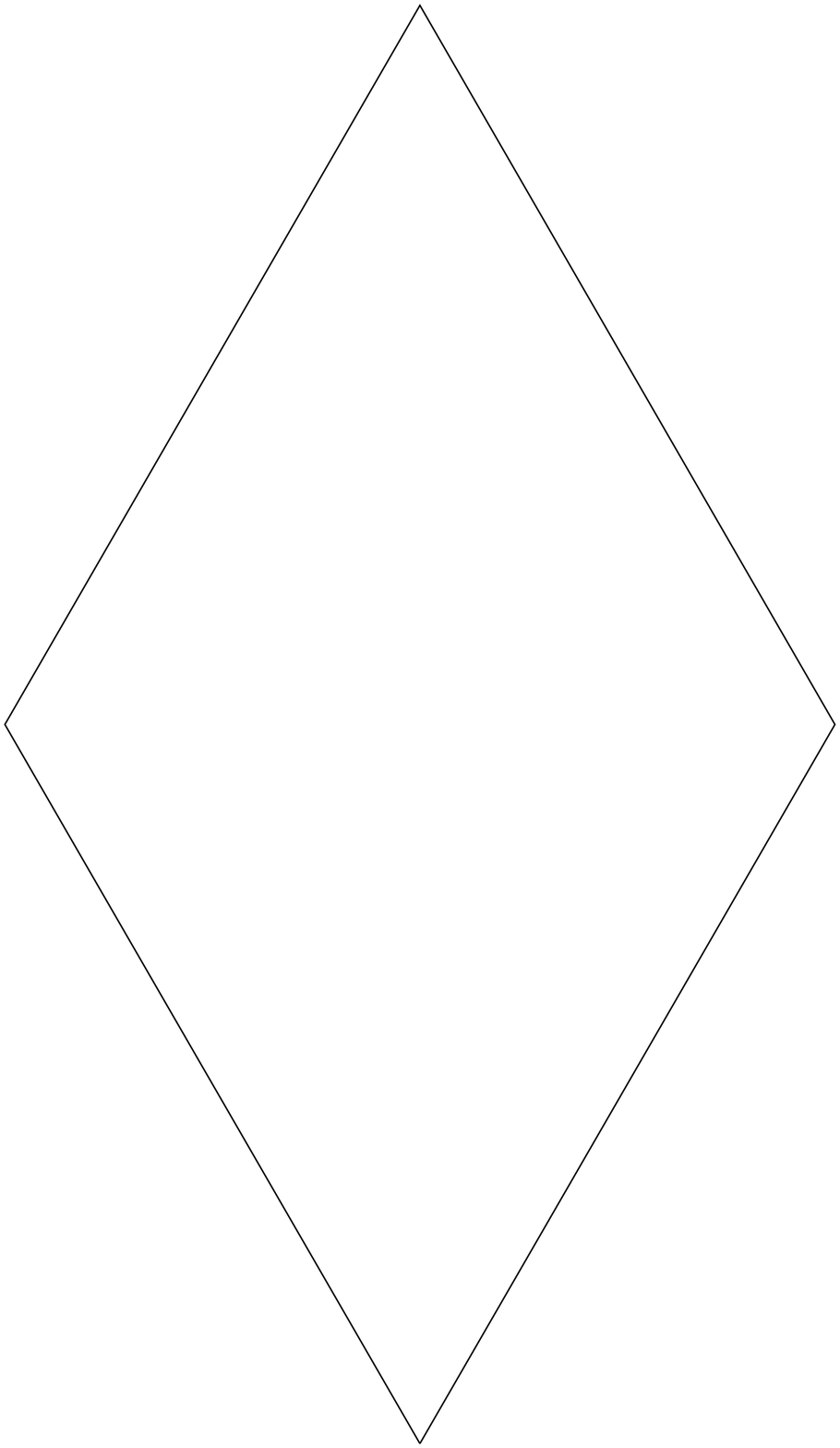}}
\put(1,0){\makebox(0,0){(a)}}
\put(2,0.5){\includegraphics[scale=1.6]{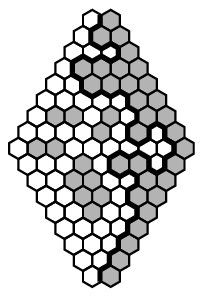}}
\put(2.9,0){\makebox(0,0){(b)}}
\put(4,0.5){\includegraphics[scale=1.6]{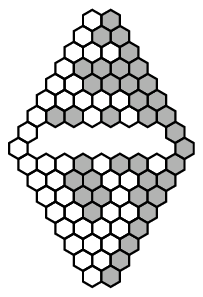}}
\put(4.9,0){\makebox(0,0){(c)}}
\put(6,0.5){\includegraphics[scale=1.6]{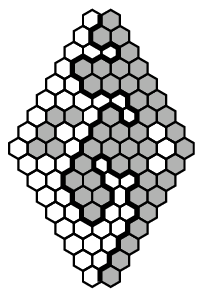}}
\put(6.9,0){\makebox(0,0){(d)}}
\end{picture}
\caption[]{\label{fig1}\small
the domain (a); percolation data and the first exploration path (b);
the percolation data near the `equator' are removed (c) and generated anew,
giving rise to the second exploration path (d).} 
\end{figure}

Being unable to prove or disprove the equivalent properties \eqref{eq1},
\eqref{eq2} for the scaling limit of percolation, I conducted an experiment as
follows. First, the domain shown on Fig.~\ref{fig1}(a) (the union of two
equilateral triangles) is filled in with percolation data. That is, each
hexagon gets black or white with probabilities $ 0.5, 0.5 $ independently of
others. However, the boundary is colored deterministically, its left part in
white, right part in black. The corresponding exploration path is
constructed. Second, the percolation data near the `equator' are removed and
generated anew. The second exploration path is constructed. Finally, the two
exploration paths are compared. This is a single trial. Many such trials are
conducted, independently of each other, for different values of parameters
(see Sect.~\ref{sect2}).

The comparison of the two exploration paths $ \ga_1, \ga_2 $ is based on the
following metric borrowed from \cite[eq.~(2)]{CM}:
\begin{equation}\label{eq3}
d(\ga_1,\ga_2) = \inf \max_{t\in[0,1]} | \ga_1(t) - \ga_2(t) | \, ;
\end{equation}
here $ | \ga_1(t) - \ga_2(t) | $ is the usual Euclidean distance between the
two points of the plane (the unit of length is such that the diameter of our
domain equals $ 2 $, that is, each side of the rhombus equals $ 2/\sqrt3 $),
and the infimum is taken over all parameterizations of the paths. Similarly to
\cite{CM}, a path is treated as an equivalence class of continuous functions $
[0,1] \to \R^2 $ modulo monotonic (increasing) re-parameterizations.

\section[]{\raggedright Experimental results}
\label{sect2}Each trial is described by two parameters $ n,k $ as follows. The number of
hexagons on the `equator' (including two boundary hexagons) equals $ n+2 $;
that is, the (discrete) rhombus consists of $ 2n+1 $ rows (of $
2,3,\dots,{n+1},\linebreak[0]
{n+2},{n+1},\dots,3,2 $ hexagons respectively). Out of
these $ 2n+1 $ rows, $ k $ rows are removed and generated anew. The set of $ k
$ rows is either symmetric w.r.t.\ the `equator' (if $ k $ is odd), or
contains one additional row above the `equator' (if $ k $ is even). For
instance, $ n=8 $ and $ k=2 $ on Fig.~\ref{fig1}(b,c,d).

Table \ref{tab1} represents $ 27 $ samples, for $ 27 $ pairs $ (n,k) $, of $
250 $ trials each. For every sample, the median of $ 250 $ values of the
distance is reported.

\begin{table}[h]
\begin{center}
\begin{tabular}{|c|ccccccc|}
\hline
\backslashbox{k}{n} & 16 & 32 & 64 & 128 & 256 & 512 & 1024 \\\hline
1 & 0.22 & 0.19 & 0.16 & 0.13 & 0.11 & 0.09 & 0.08 \\
2 & & 0.24 & 0.18 & 0.19 & 0.15 & 0.11 & 0.10 \\
4 & & & 0.27 & 0.24 & 0.19 & 0.15 & 0.12 \\
8 & & & & 0.28 & 0.23 & 0.19 & 0.16 \\
16 & & & & & 0.29 & 0.22 & 0.18 \\
32 & & & & & & 0.31 & 0.22 \\\hline
\end{tabular}
\end{center}
\caption{\label{tab1}\small
Median values of the distance between the two exploration paths.}
\end{table}

The distance is basically defined by \eqref{eq3}. However, for speeding up the
computation, I use an approximation to \eqref{eq3} (described in
Sect.~\ref{sect4}), which leads to a systematic error bounded by $ \pm 0.03
$. Thus, each median in Table \ref{tab1} has a systematic error bounded by $
\pm 0.03 $, and a sampling error, the mean square deviation being about $ 7 \%
$ of the median.

The results suggest that the typical distance tends to $ 0 $ as the strip
width $ \eps = \frac k n $ tends to $ 0 $. Maybe we observe a power law $
\eps^\al $ for some exponent $ \al $ not far from $ 1/3 $.

\section[]{\raggedright Examples}
\label{sect3}
\begin{figure}[p]
\setlength{\unitlength}{2.3cm}
\begin{picture}(7,5)
\put(0,2.3){\includegraphics[scale=0.2]{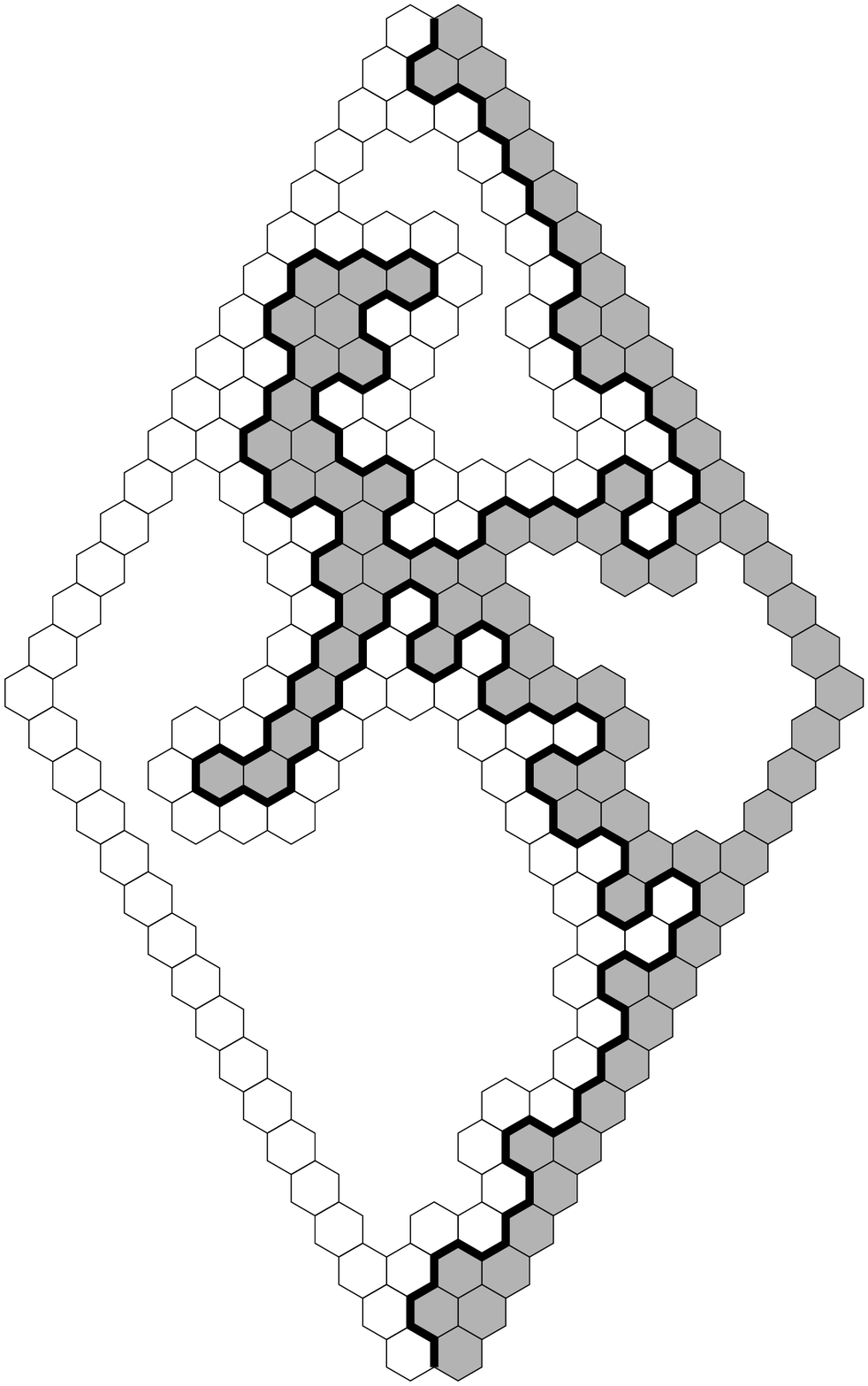}}
\put(0.75,2){\makebox(0,0){(a)}}
\put(2,2.3){\includegraphics[scale=0.2]{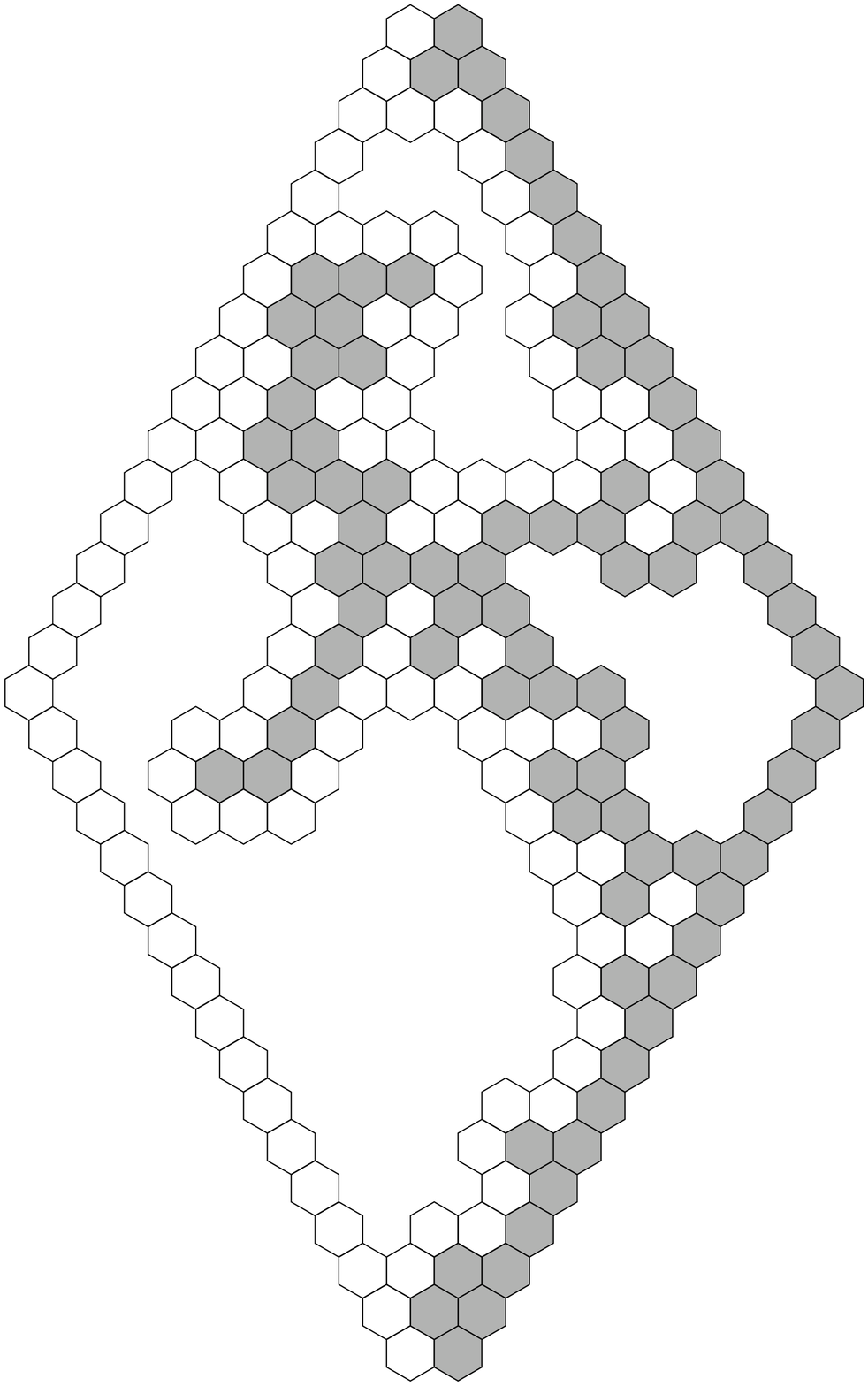}}
\put(2.75,2){\makebox(0,0){(b)}}
\put(4,2.3){\includegraphics[scale=0.2]{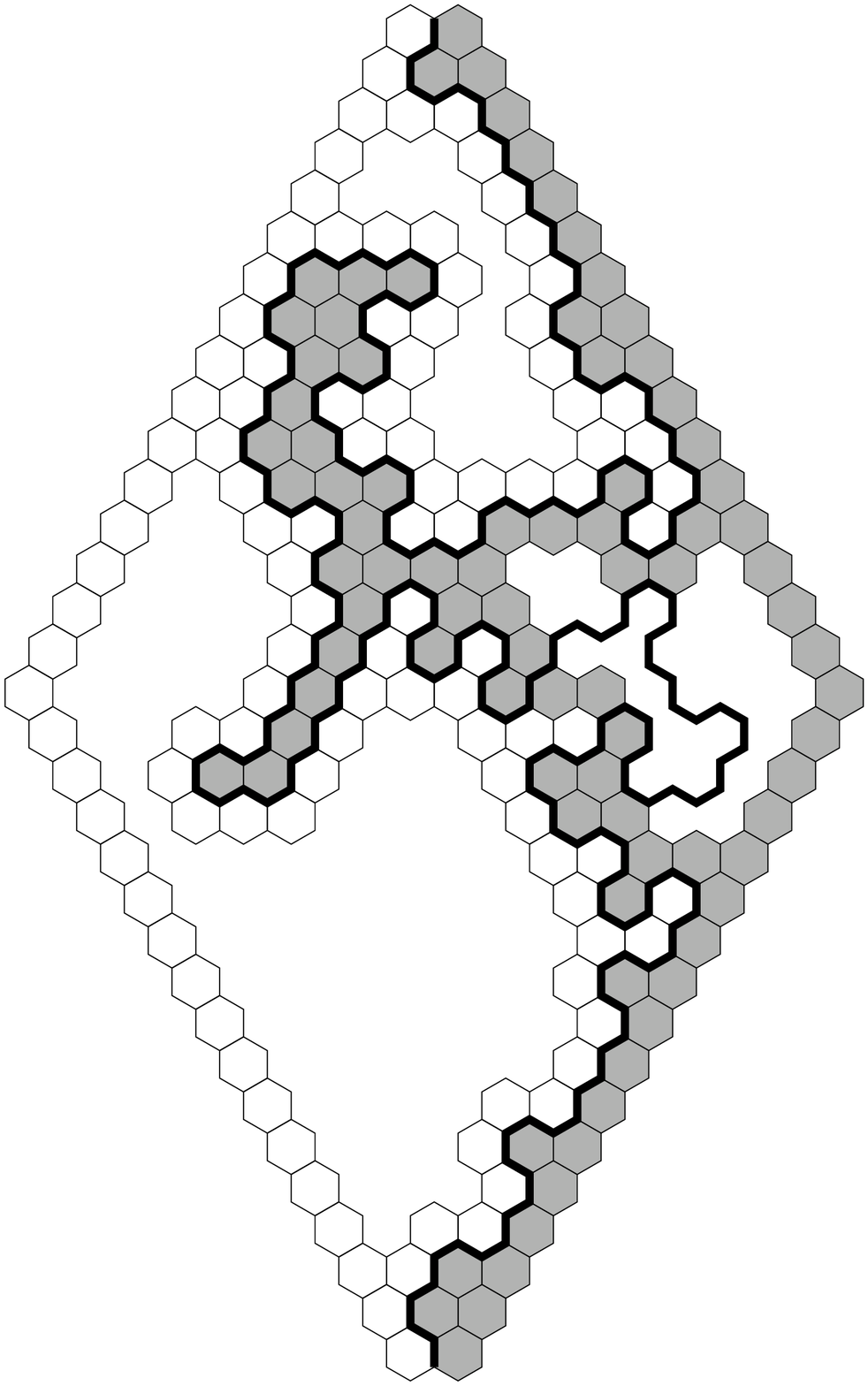}}
\put(4.75,2){\makebox(0,0){(c)}}
\put(1.5,0.3){\includegraphics[scale=0.35,bb=-227 248 254 488,clip]{pic2c.eps}}
\put(2.75,0){\makebox(0,0){(d)}}
\end{picture}
\caption[]{\label{fig2}\small
The median trial for $ n=16 $, $ k=1 $. The first exploration path (a) and
the hexagons relevant to it (a,b); the same hexagons together with the second
exploration path (c,d).}
\end{figure}

The median ($ =0.22 $) of the first sample ($ n=16, k=1 $) resulted from a
definite trial (no.~21, in fact), --- call it the \emph{median trial,} ---
presented on Fig.~\ref{fig2}. Only the first exploration path is shown on
Fig.~\ref{fig2}(a) together with the hexagons relevant to it. On
Fig.~\ref{fig2}(b) the path is not shown explicitly, but you can restore it
easily. On Fig.~\ref{fig2}(c) the \emph{second} exploration path is shown
together with the hexagons relevant to the \emph{first} path. Hexagons
relevant to the second path are not shown explicitly, but you can restore them
easily. Fig.~\ref{fig2}(d) presents the most interesting, `equatorial' part of Fig.~\ref{fig2}(c).

\begin{figure}[p]
\setlength{\unitlength}{2.3cm}
\begin{picture}(7,2)
\put(0,0.3){\includegraphics[scale=0.4,bb=-227 238 254 498,clip]{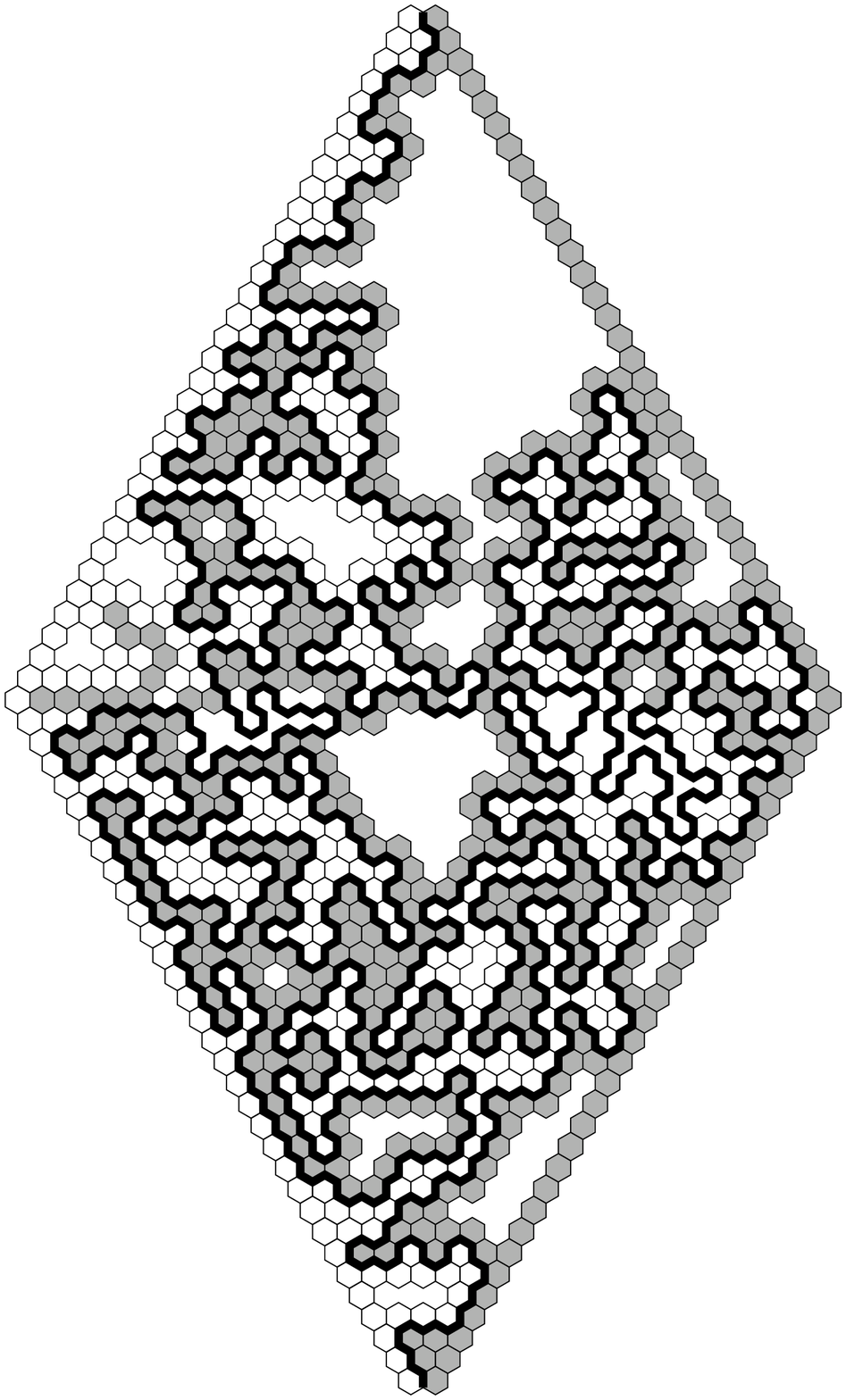}}
\put(1.5,0){\makebox(0,0){$ k=1 $}}
\put(3,0.3){\includegraphics[scale=0.4,bb=-227 238 254 498,clip]{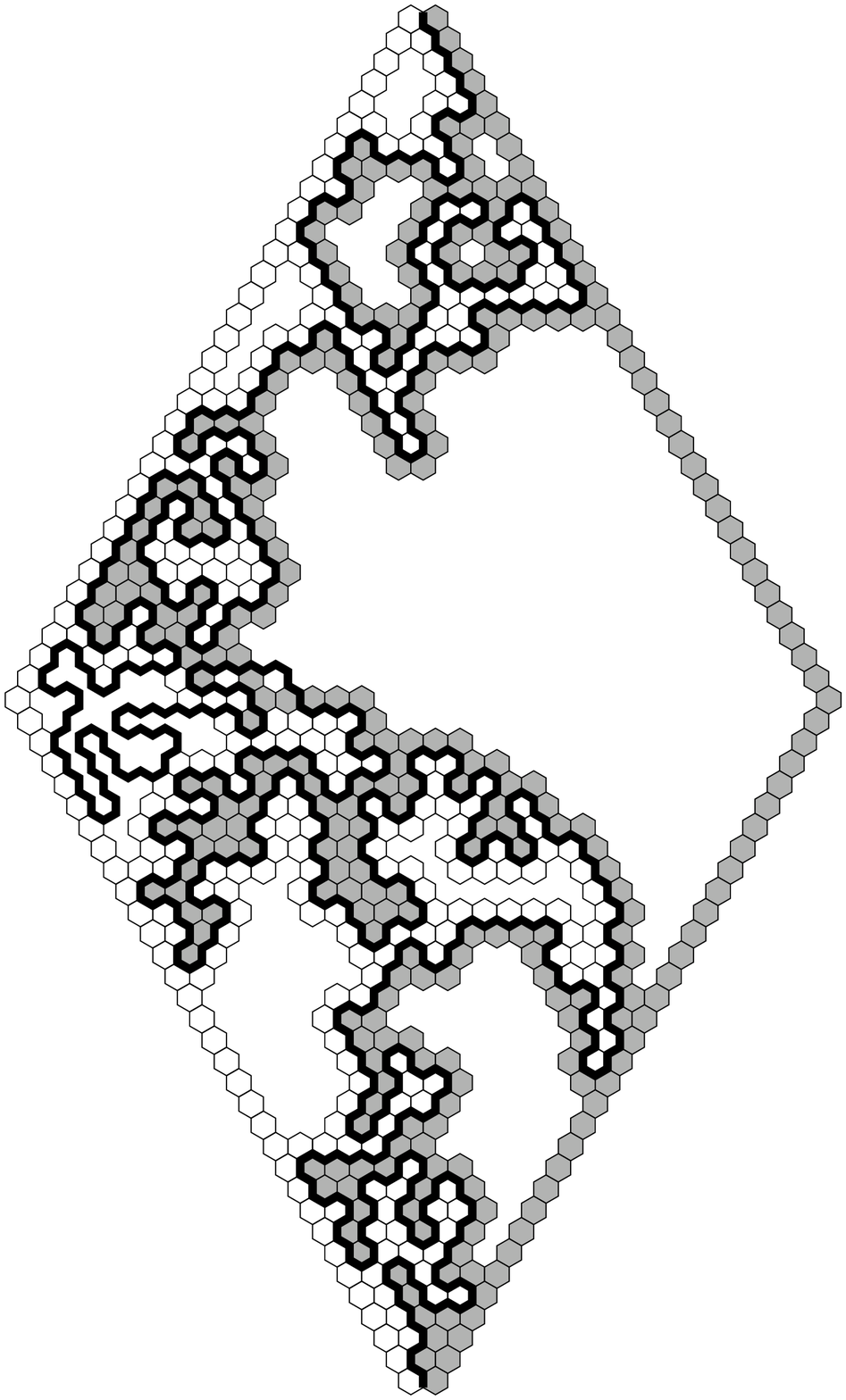}}
\put(4.5,0){\makebox(0,0){$ k=2 $}}
\end{picture}
\caption[]{\label{fig3}\small
Median trials for $ n=32 $.}
\end{figure}

\begin{figure}[p]
\setlength{\unitlength}{2.3cm}
\begin{picture}(7,8.9)
\put(0,7.2){\includegraphics[scale=0.8,bb=-227 288 254 448,clip]{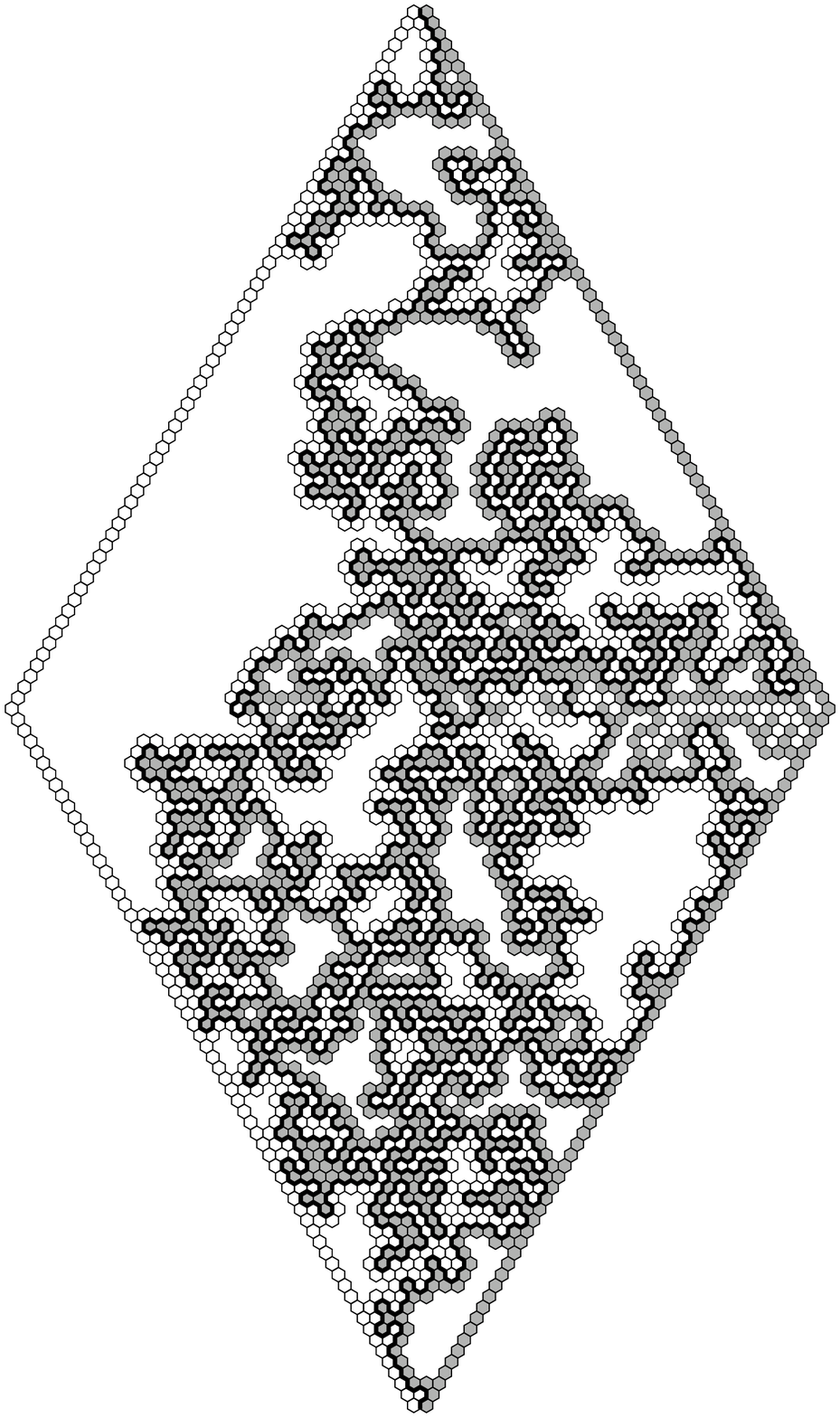}}
\put(3,7.0){\makebox(0,0){$ k=1 $}}
\put(0,3.9){\includegraphics[scale=0.8,bb=-227 258 254 478,clip]{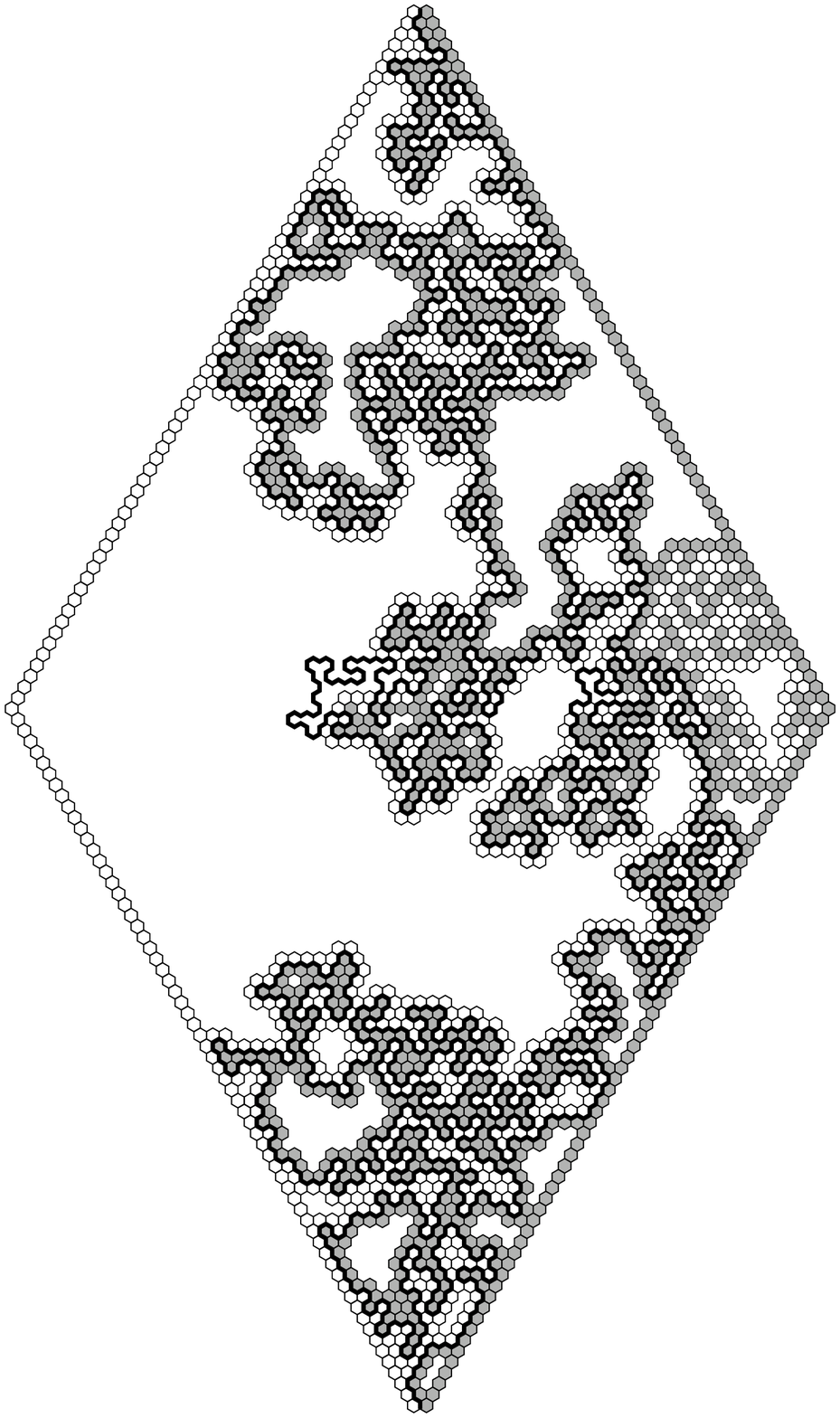}}
\put(3,4.0){\makebox(0,0){$ k=2 $}}
\put(0,0.2){\includegraphics[scale=0.8,bb=-227 228 254 508,clip]{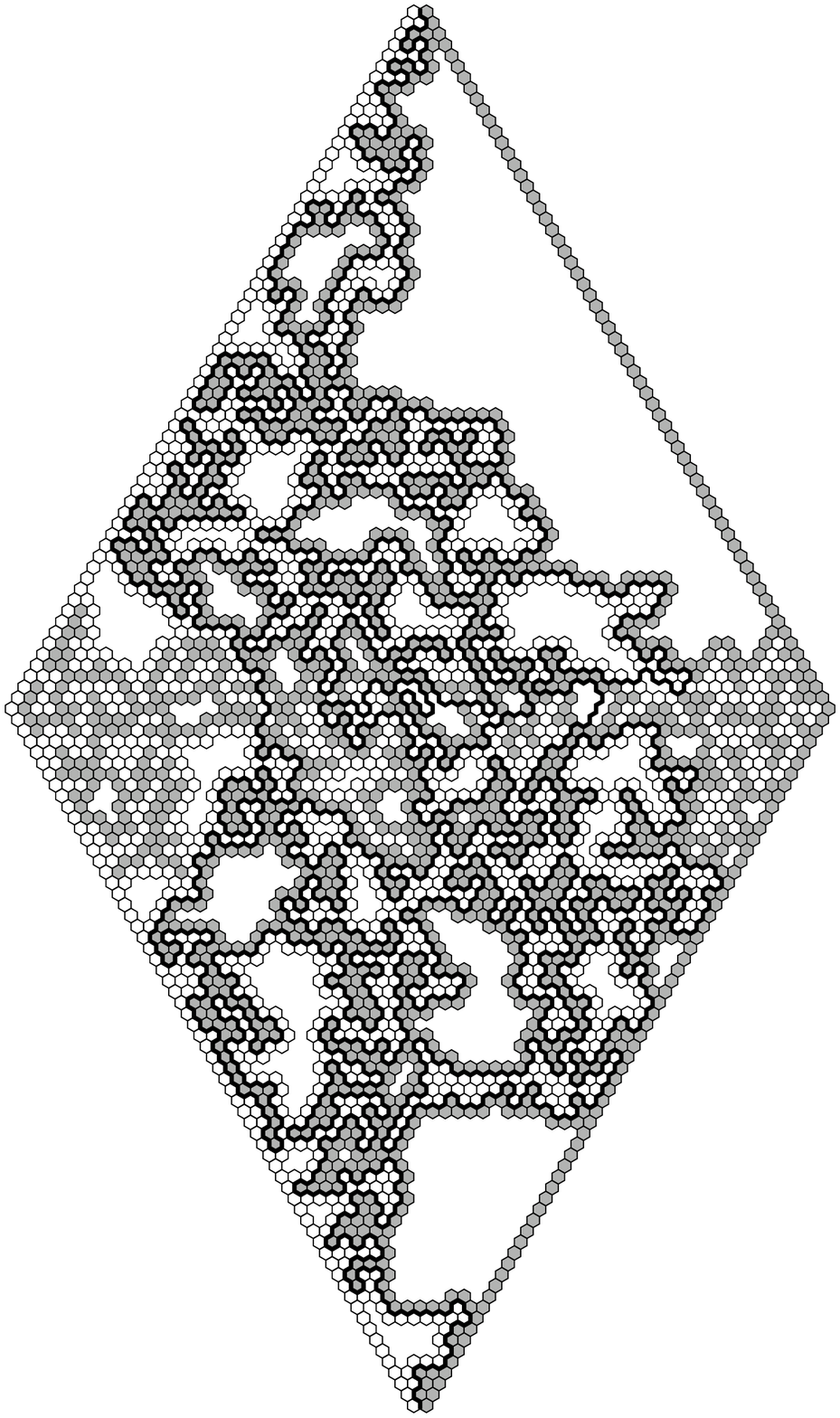}}
\put(3,0){\makebox(0,0){$ k=4 $}}
\end{picture}
\caption[]{\label{fig4}\small
Median trials for $ n=64 $.}
\end{figure}

More examples follow (Fig.~\ref{fig3}, \ref{fig4}), in the same
format as Fig.~\ref{fig2}(d).

\section[]{\raggedright Technicalities}
\label{sect4}First, about the metric \eqref{eq3} and its discrete approximation
used. Treating the exploration path as a finite sequence of points, ---
vertices of the polygonal line (rather than the line itself), I define a
distance $ d(a,b) $ between two such sequences $ a = (a_0,\dots,a_M) $ and $ b
= (b_0,\dots,b_N) $ as follows:
\[
d(a,b) = \min_{(\mu,\nu)} \max_{k=0,\dots,M+N} | a_{\mu_k} - b_{\nu_k} | \, ;
\]
here $ | a_{\mu_k} - b_{\nu_k} | $ is the Euclidean distance between the two
points, and $ (\mu,\nu) $ runs over all pairs of sequences $
\mu_0,\dots,\mu_{M+N} $ and $ \nu_0,\dots,\nu_{M+N} $ of integers satisfying
\begin{gather*}
0 = \mu_0 \le \dots \le \mu_{M+N} = M \, , \quad
0 = \nu_0 \le \dots \le \nu_{M+N} = N \, , \\
\mu_k + \nu_k = k \quad \text{for } k=0,\dots,M+N \, .
\end{gather*}
Clearly it is an approximation to \eqref{eq3}. (Not a metric, since $ d(a,a)
\ne 0 $, which is harmless.) Dynamic programming helps us to compute the
minimum over $ (\mu,\nu) $. To this end, introduce
\[
d_{k,i} = \min_{(\mu,\nu):\mu_k=i} \max_{l=0,\dots,k} | a_{\mu_l} - b_{\nu_l}
|
\]
for $ i = 0,\dots,\min(k,M) $ and consider the array $ D_k = \( d_{k,i}
\)_{i=0,\dots,\min(k,M)} $. Having $ D_k $ we can compute $ D_{k+1} $
easily. Thus, we compute $ D_0 $, then $ D_1 $, and so on, up to $ D_{M+N} $;
the latter gives us $ d(a,b) $.

However, $ M+N $ is typically as large as $ 10^5 $ (for $ n=1024 $, my worst
case). For speeding up the computation, I replace the long sequence $ a =
(a_0,\dots,a_M) $ with a subsequence $ (a_{k_1},\dots,a_{k_p}) $ such that $
|a_{k_i} - a_l | \le \eps $ for $ k_i < l < k_{i+1} $. (The same is done for
the other sequence $ b $.) The parameter $ \eps $ was chosen to be $ 0.03 $.

Now, other technicalities.
Random colors (of the hexagons) are produced by a Wichmann-Hill generator of
pseudo-random numbers, initialized (before each trial) with a one-to-one
function of $ \log_2 n $, $ \log_2 k $ and the number of the trial (within the
sample).

The program, in the `Python' programming language, is available on my Web
site.

The computer: a personal computer with matherboard Intel D815EGEW, processor
Intel Pentium III (933 MHz), cache RAM 256 Kb and SDRAM 256 Mb (133 MHz).

The operating system: SuSE Linux 9.0.

Python: version 2.3+.

One trial takes typically from 0.5 min to 1 min for $ n=1024 $, and 10--15 sec
for $ n=512 $.

\section*{Conclusion}

The results (Table \ref{tab1}) suggest that the scaling limit of percolation
is \emph{strip stable,} and therefore may be treated as a black noise.

One could object that the strong concentrated perturbation is not enough; a
distributed `microscopic' perturbation should be added, such that exploration
paths remain macroscopically close to the same \SLE\ curves, but are changed
microscopically. That is correct, but I guess that such `micro' perturbations
may be absorbed somehow by the strong perturbation in the infinitesimal strip.

\bigskip
\filbreak
{
\small
\begin{sc}
\parindent=0pt\baselineskip=12pt
\parbox{4in}{
Boris Tsirelson\\
School of Mathematics\\
Tel Aviv University\\
Tel Aviv 69978, Israel
\smallskip
\par\quad\href{mailto:tsirel@post.tau.ac.il}{\tt
 mailto:tsirel@post.tau.ac.il}
\par\quad\href{http://www.tau.ac.il/~tsirel/}{\tt
 http://www.tau.ac.il/\textasciitilde tsirel/}
}

\end{sc}
}
\filbreak

\end{document}